\documentclass[12pt]{article}
\usepackage{amsmath, amssymb, amsthm}
\usepackage{graphicx}

\newtheorem{theorem}{Theorem}%[section]

\newtheorem{algorithm}[theorem]{Algorithm}

\newtheorem*{general Gromov'}{Corollary \ref{general Gromov}$'$}

\def \proof {\noindent {\bf Proof.}\ \ }

\def \endproof {{\mbox{}\nolinebreak\hfill\rule{2mm}{2mm}\par\medbreak}}

\def\liminf{{\operatorname{liminf}}}
\def\limsup{{\operatorname{limsup}}}

\def \R {\mathbb{R}}
\def \C {\mathbb{C}}
\def \Cst {\mathbb{C}}

\def \E {\mathbb{E}}

\def \e {\varepsilon}
\def \eps {\varepsilon}

\def \k {\kappa}

\def \< {\langle}
\def \> {\rangle}

\def \lin {{\rm span}}

\def\eps{{\varepsilon}}
\def\Ae{A^{(e)}}
\def\be{b^{(e)}}
\def\ord{{\mathcal O}}

\begin{document}
%\title {A randomized solver for linear systems with exponential convergence}
\title {A randomized Kaczmarz algorithm with exponential convergence}
\author{Thomas Strohmer and Roman Vershynin\footnote{T.S.\ was supported
by NSF DMS grant 0511461. R.V. was supported by the Alfred P.\ Sloan 
Foundation and by NSF DMS grant 0401032.}\\
Department of Mathematics, University of California \\
Davis, CA 95616-8633, USA.\\
strohmer@math.ucdavis.edu, vershynin@math.ucdavis.edu}

\date{}
\maketitle

\begin{abstract}
The Kaczmarz method for solving linear systems of equations is
an iterative algorithm that has found many applications ranging from
computer tomography to digital signal processing. Despite the popularity
of this method, useful theoretical estimates for its rate of convergence
are still scarce. We introduce a randomized version of the
Kaczmarz method for consistent, overdetermined linear systems and we prove
that it converges with expected exponential rate. 
Furthermore, this is the first solver whose {\em rate does not depend on 
the number of equations} in the system.  
The solver does not even need to know the whole system, but only a small 
random part of it.
It thus outperforms all previously known methods on general 
extremely overdetermined systems.
Even for moderately overdetermined systems, numerical simulations as
well as theoretical analysis reveal that our algorithm can converge faster 
than the celebrated conjugate gradient algorithm.  
Furthermore, our theory and numerical simulations confirm a prediction
of Feichtinger et al.\ in the context of reconstructing bandlimited
functions from nonuniform sampling.
\end{abstract}

\section{Introduction and state of the art} \label{intro}
%---------------------------------------------------------------
We study a consistent linear system of equations
\begin{equation}
\label{system}
A x = b,
\end{equation}
where $A$ is a full rank $m \times n$ matrix with $m \ge n$,
and $b\in \Cst^m$. One of the most popular solvers for such overdetermined systems
is {\em Kaczmarz's method}~\cite{Kac37}, which is a form of
alternating projection method. This method is also known under the name
{\em Algebraic Reconstruction Technique} (ART) in
computer tomography~\cite{Her80,Nat86}, and in fact, it was implemented in
the very first medical scanner~\cite{Hou73}. It can also be considered as a
special case of the POCS (Projection onto Convex Sets) method, which is a
prominent tool in signal and image processing~\cite{SS87,CFM92}.

%Let the rows of $A$ be represented by the vectors
%$a_1,\dots,a_m \in \Cst^n$, and let $b=(b_1,\dots,b_m)^T$. 
We denote the rows of $A$ by $a_1^{\ast},\dots,a_m^{\ast}$
and let $b=(b_1,\dots,b_m)^T$. 
The classical scheme of Kaczmarz's method sweeps
through the rows of $A$ in a cyclic manner, projecting in each substep
the last iterate orthogonally onto the solution hyperplane
of $\langle a_i,x\rangle = b_i$ and taking this as the next iterate. Given
some initial approximation $x_0$, the algorithm takes the form
\begin{equation}
x_{k+1} = x_k +\frac{b_i - \langle a_i,x_k\rangle}{\|a_i\|_2^2}a_i,
\label{kacz}
\end{equation}
where $i=k \mod m+1$ and $\|\cdot\|$ denotes the Euclidean norm 
in $\C^n$.
% In this paper we will restrict ourselves
%to the case where the relaxation parameter $\lambda_i = 1$ for all $i$.
Note that we refer to one projection as one iteration, thus
one sweep in~\eqref{kacz} through all $m$ rows of $A$ consists of $m$
iterations. 

While conditions for convergence of this method are readily established,
useful theoretical estimates of the {\em rate of convergence} of
the Kaczmarz method (or more generally of the alternating projection method
for linear subspaces) are difficult to obtain, at least for $m > 2$.
Known estimates for the rate of convergence are based on quantities 
of the matrix $A$ that are hard to compute and 
difficult to compare with convergence estimates of other iterative methods
(see e.g.~\cite{Deu85,DH97,Gal05} and the references therein).

What numerical analysts would like to have is estimates of the convergence
rate in terms of a {\em condition number} of $A$.
No such estimates have been known prior to this work.
The difficulty stems from the fact that the rate of convergence 
of~\eqref{kacz} depends strongly on the
{\em order} of the equations in~\eqref{system}, while 
condition numbers do not depend on the order of the rows of a matrix.

It has been observed several times in the literature that using the rows
of $A$ in Kaczmarz's method in random order, rather than in their given order,
can greatly improve the rate of convergence, see e.g.~\cite{Nat86,CFM92,HM93}.
While this randomized Kaczmarz method is thus quite appealing for
applications, no guarantees of its rate of convergence have been known.

In this paper, we propose the first randomized Kaczmarz method
with {\em exponential expected rate of convergence}, 
cf.~Section~\ref{s:random}. 
Furthermore, this rate depends only on the scaled condition number
of $A$ and {\em not on the number of equations} $m$ in the system.
The solver does not even need to know the whole system, but only a small
random part of it. Thus our solver outperforms all previously known methods on general extremely overdetermined systems.

We analyze the optimality of the proposed algorithm as well as of the 
derived estimate, cf.~Section~\ref{s:optim}.
Section~\ref{s:numeric} contains various numerical simulations. In one
set of experiments we apply the randomized Kaczmarz method to the 
reconstruction of bandlimited functions from non-uniformly spaced samples.
In another set of numerical simulations, accompanied by theoretical
analysis, we demonstrate that even for moderately overdetermined systems, 
the randomized Kaczmarz method can outperform the celebrated conjugate 
gradient algorithm.

\medskip
{\bf Condition numbers.}
For a matrix $A$, its spectral norm is denoted by $\|A\|_2$,
and its Frobenius norm by $\|A\|_F$.
Thus the spectral norm is the
largest singular value of $A$, and the Frobenius norm is the 
square root of the sum of the squares of all singular values of $A$.

The left inverse of $A$ (which we always assume to exist) is 
denoted by $A^{-1}$. Thus $\|A^{-1}\|_2$ is the smallest 
constant $M$ such that the inequality 
$\|Ax\|_2 \ge \frac{1}{M} \|x\|_2$ holds for all vectors $x$.

The usual condition number of $A$ is 
$$
k(A):=\|A\|_2 \|A^{-1}\|_2.
$$
A related version is the scaled condition number 
introduced by Demmel \cite{D}:
$$
\k(A):=\|A\|_F \|A^{-1}\|_2.
$$
One easily checks that 
\begin{equation}				\label{condition numbers}
  1 \le \frac{\kappa(A)}{\sqrt{n}} \le k(A).
\end{equation}
Estimates on the condition numbers of some typical (i.e. random or
Toeplitz-type ) matrices are known from a big body of literature, 
see \cite{BG,D, Ede88, Ede92,  ES, RV1, RV2,ST,TV} and the references therein.

\section{Randomized Kaczmarz algorithm and its rate of convergence}
\label{s:random}

It has been observed in numerical simulations~\cite{Nat86,CFM92,HM93} that 
the convergence rate of the Kaczmarz method can be significantly improved 
when the algorithm~\eqref{kacz} sweeps through the rows of $A$ in a random 
manner, rather than sequentially in the given order. In fact, the improvement 
in convergence can be quite dramatic. Here we propose a specific
version of this randomized Kaczmarz method, which chooses each row of $A$ 
with probability proportional to its relevance -- more precisely, 
proportional to the square of its Euclidean norm.
This method of sampling from a matrix was proposed in \cite{FKV98}
in the context of computing a low-rank approximation of $A$,
see also \cite{RV} for subsequent work and references.
Our algorithm thus takes the following form:
\begin{algorithm}[Random Kaczmarz algorithm]
\label{alg}
Let $Ax=b$ be a linear system of equations as in~\eqref{system} and
let $x_0$ be arbitrary initial approximation to the solution
of~\eqref{system}. For $k=0,1,\dots$ compute
\begin{equation}
\label{kaczrand}
x_{k+1} = x_k +\frac{b_{r(i)} -\langle a_{r(i)},x_k\rangle}
                    {\|a_{r(i)}\|_2^2} a_{r(i)},
\end{equation}
where $r(i)$ is chosen from the set $\{1,2,\dots,m\}$
at random, with probability proportional to $\|a_{r(i)}\|_2^2$.
\end{algorithm}

\if 0
\begin{algorithm}[Random Kaczmarz algorithm]
\label{alg}
Let $Ax=b$ be a linear system of equations as in~\eqref{system}.
\begin{itemize}
\item Let $x_0$ be arbitrary initial approximation to the solution
of~\eqref{system}.
\item At step $k = 1, 2, \ldots$, choose one equation from the
system \eqref{system} at random, with probability proportional to the
square of the Euclidean length of the corresponding row of the matrix $A$.
Let $x_k$ be the orthogonal projection of the running
approximation $x_{k-1}$ onto the solution space of that equation.
\item Output $x_k$ when the desired precision is obtained.
\end{itemize}
\end{algorithm}
\fi

Our main result states that $x_k$ converges exponentially
fast to the solution of~\eqref{system}, and the rate of convergence
depends only on the scaled condition number $\k(A)$.

\begin{theorem}
\label{solution}
Let $x$ be the solution of \eqref{system}. Then Algorithm~\ref{alg}
converges to $x$ in expectation, with the average error
\begin{equation}
\label{mainest}
\E \|x_k - x\|_2^2 \le \bigl(1 - \k(A)^{-2} \bigr)^k \cdot \|x_0 - x\|_2^2.
\end{equation}
\end{theorem}

\begin{proof}
There holds
\begin{equation}
\sum_{j=1}^m |\< z, a_j \> |^2 \ge  \frac{\|z\|_2^2}{\|A^{-1}\|_2^2}
 \qquad \text{for all $z \in \C^n$.}        \label{rows of A}
\end{equation}
Using the fact that $\|A\|_F^2=\sum_{j=1}^m \|a_j\|_2^2$ we can
write~\eqref{rows of A} as
\begin{equation}                \label{normalized rows}
\sum_{j=1}^m \frac{\|a_j\|_2^2}{\|A\|_F^2} \;
  \Bigl| \Big\langle z, \frac{a_j}{\|a_j\|_2} \Big\rangle \Big|^2
\ge \k(A)^{-2} \|z\|^2
 \ \ \ \text{for all $z \in \C^n$}.
\end{equation}
The main point of the proof is to view the left hand side
in~\eqref{normalized rows} as an expectation of some random variable.
Namely, recall that the solution space of the $j$-th equation of \eqref{system}
is the hyperplane $\{y : \< y, a_j \> = b_j \}$, whose normal is
$\frac{a_j}{\|a_j\|_2}$.
Define a random vector $Z$ whose values are the normals to all
the equations of \eqref{system}, with probabilities as in our algorithm:
\begin{equation}                    \label{Z}
Z  =  \frac{a_j}{\|a_j\|_2}
\ \ \ \text{with probability} \ \ \
\frac{\|a_j\|_2^2}{\|A\|_F^2}, \ \ \
j = 1, \ldots, m.
\end{equation}
Then \eqref{normalized rows} says that
\begin{equation}                    \label{expectation}
\E | \< z, Z\> |^2  \ge \k(A)^{-2} \|z\|_2^2
 \ \ \ \text{for all $z \in \C^n$.}
\end{equation}
The orthogonal projection $P$ onto the solution space of a random
equation of \eqref{system} is given by
$P z = z - \< z - x, Z \> \, Z$.

Now we are ready to analyze our algorithm. We want to show that
the error $\|x_k - x\|_2^2$ reduces at each step in average
(conditioned on the previous steps) by at least the factor
of $(1 - \k(A)^{-2})$.
The next approximation $x_k$ is computed from $x_{k-1}$ as
$x_k = P_k x_{k-1}$, where $P_1, P_2, \ldots$ are independent realizations
of the random projection $P$.
The vector $x_{k-1} - x_k$ is in the kernel of $P_k$. It is
orthogonal to the solution space of the equation onto which $P_k$
projects, which contains the vector $x_k - x$ (recall that $x$ is the
solution to all equations). The orthogonality of these two vectors then yields
$$
\|x_k - x\|_2^2 = \|x_{k-1} - x\|_2^2 - \|x_{k-1} - x_k\|_2^2.
$$
To complete the proof, we have to bound $\|x_{k-1} - x_k\|_2^2$ from below.
By the definition of $x_k$, we have
$$
\|x_{k-1} - x_k\|_2  =  \< x_{k-1} - x, Z_k \>
$$
where $Z_1, Z_2, \ldots$ are independent realizations of the random vector $Z$.
Thus
$$
\|x_k - x\|_2^2
\le \Big( 1 - \Big| \Big\langle \frac{x_{k-1} - x}{\|x_{k-1} - x\|_2},
                  Z_k
                  \Big\rangle \Big|^2 \Big)
  \;\|x_{k-1} - x\|_2^2.
$$
Now we take the expectation of both sides conditional upon
the choice of the random vectors $Z_1, \ldots, Z_{k-1}$
(hence we fix the choice of the random projections $P_1, \ldots, P_{k-1}$ and
thus the random vectors $x_1, \ldots, x_{k-1}$, and we average over the 
random vector $Z_k$). Then
$$
\E_{\{Z_1,\ldots,Z_{k-1}\}}  \|x_k - x\|_2^2
\le \Big( 1 - \E_{\{Z_1,\ldots,Z_{k-1}\}}
                     \Big| \Big\langle \frac{x_{k-1} - x}{\|x_{k-1} - x\|_2},
                     Z_k
                     \Big\rangle \Big|^2 \Big)
  \; \|x_{k-1} - x\|_2^2.
$$
By \eqref{expectation} and the independence,
$$
\E_{\{Z_1,\ldots,Z_{k-1}\}}  \|x_k - x\|_2^2
\le   \big( 1 - \k(A)^{-2} \big) \; \|x_{k-1} - x\|_2^2.
$$
Taking the full expectation of both sides, we conclude that 
$$
\E \|x_k - x\|_2^2
\le   \big( 1 - \k(A)^{-2} \big) \; \E \|x_{k-1} - x\|_2^2.
$$
By induction, we complete the proof.
\end{proof}

\subsection{Quadratic time}
Theorem~\ref{solution} yields a simple bound on the expected 
computational complexity of the randomized Kaczmarz Algorithm~\ref{alg} 
to compute the solution within error $\eps$, i.e. 
\begin{equation}			\label{within eps}
  \E \|x_k - x\|_2^2 \le \e^2 \|x_0 - x\|_2^2.
\end{equation}
The expected number of iterations (projections) $k_\e$ to achieve an accuracy 
$\eps$ is
\begin{equation}				\label{k-eps}
  \E \, k_\e \le \frac{2 \log \e}{\log (1 - \k(A)^{-2})} 
  \approx 2 \k(A)^2 \log \frac{1}{\e},
\end{equation}
where $f(n) \sim g(n)$ means $f(n)/g(n) \to 1$ as $n \to \infty$.
(Note that $\k(A)^2 \ge n$ by \eqref{condition numbers}, so the 
approximation in \eqref{k-eps} becomes tight as the number of equation
$n$ grows).

Since each projection can be computed in $O(n)$ time,
and $\kappa(A)^2 = O(n)$ by \eqref{condition numbers},
{\em the algorithm takes $O(n^2)$ operations to converge to the solution}.
This should be compared to the Gaussian elimination, which takes $O(mn^2)$
time. Strassen's algorithm and its improvements reduce the exponent in
Gaussian elimination, but these algorithms are, as of now, of no practical use.

%Note that our algorithm does not compute the inverse of $A$,
%only its action on one arbitrary vector $b$.
Of course, we have to know the (approximate) Euclidean lengths of the rows
of $A$ before we start iterating; computing them takes $O(nm)$ time.
But the lengths of the rows may in many cases be known a priori.
For example, all of them may be equal to one (as is the case for
Vandermonde matrices arising in trigonometric approximation) or tightly
concentrated around a constant value (as is the case for random
matrices).

The number $m$ of the equations is essentially irrelevant for our algorithm. 
The algorithm does not even need to know the whole matrix, 
but only $O(n)$ random rows.
Such Monte-Carlo methods have been successfully developed for many
problems, even with precisely the same model of selecting a random submatrix
of $A$ (proportional to the squares of the lengths of the rows),
see \cite{FKV98} for the original discovery and \cite{RV} for subsequent
work and references.

\section{Optimality}
\label{s:optim}
%-------------------------------------------------

We discuss conditions under which our algorithm is optimal in
a certain sense, as well as the optimality of the estimate
on the expected rate of convergence.

\subsection{General lower estimate} \label{ss:lowerest}

For any system of linear equations, our estimate can not
be improved beyond a constant factor, as shown by the following
theorem.
\begin{theorem}			\label{general lower estimate}
Consider the linear system of equations~\eqref{system} and let $x$
be its solution. Then there exists an initial approximation $x_0$
such that
\begin{equation}
\label{optimalR}
\E \|x_k - x\|_2^2 
\ge \big( 1 - 2k / \k(A)^2 \big) \cdot \|x_0 - x\|_2^2
\end{equation}
for all $k=1,2,\ldots$
\end{theorem}

\medskip

\proof
For this proof we can assume without loss of generality that the
system~\eqref{system} is homogeneous: $Ax = 0$. Let $x_0$ be a vector which
realizes $\kappa(A)$, that is
$\kappa(A) = \|A\|_F \|A^{-1} x_0\|_2$ and $\|x_0\|_2 = 1$.
As in the proof of Theorem~\ref{solution}, we define the random normal $Z$
associated with the rows of $A$ by \eqref{Z}.
%by $a^{\ast}_1, \ldots, a^{\ast}_m$ and
Similar to \eqref{expectation}, we have
$\E | \< x_0, Z\> |^2  = \k(A)^{-2}$.
We thus see $\lin(x_0)$ as an ``exceptional'' direction, so we shall
decompose $\R^n = \lin(x_0) \oplus (x_0)^\perp$, writing every vector
$x \in \R^n$ as
$$
x = x' \cdot x_0 + x'',
\ \ \ \text{where} \ \ \
x' \in \R, \ \ x'' \in (x_0)^\perp.
$$
In particular,
\begin{equation}                    \label{expectation Z'}
\E |Z'|^2 = \kappa(A)^{-2}.
\end{equation}

We shall first analyze the effect of one random projection in our
algorithm. To this end, let $x \in \R^n$, $\|x\|_2 \le 1$,
and let $z \in \R^n$, $\|z\|_2 = 1$.
(Later, $x$ will be the running approximation $x_{k-1}$, and $z$
will be the random normal $Z$).
The projection of $x$ onto the hyperplane whose normal is $z$ equals
$$
x_1 = x - \< x, z\> z.
$$
Since
\begin{equation}                    \label{scalar product decomposed}
\< x, z \>  =  x' z' + \< x'', z'' \> ,
\end{equation}
we have
\begin{equation}                    \label{prime}
|x'_1 -x'|
= |\< x, z \> z'|
\le |x'||z'|^2 + |\< x'', z'' \> z'|
\le |z'|^2 + |\< x'', z''\> z'|
\end{equation}
because $|x'| \le \|x\|_2 \le 1$.
Next,
\begin{align*}
\|x''_1\|^2 - \|x''\|^2
&= \|x'' - \< x, z\> z'' \|_2^2 - \|x''\|_2^2 \\
&= -2 \< x, z \> \< x'', z'' \> + \< x, z\> ^2 \|z''\|_2^2 
\le  -2 \< x, z \> \< x'', z'' \> + \< x, z\> ^2
\end{align*}
because $\|z''\|_2 \le \|z\|_2 = 1$.
Using \eqref{scalar product decomposed}, we decompose
$\< x, z \> $ as $a + b$, where $a = x'z'$ and $b = \< x'', z'' \> $
and use the identity $-2(a+b)b + (a+b)^2 = a^2 - b^2$
to conclude that
\begin{equation}                    \label{double prime}
\|x''_1\|_2^2 - \|x''\|_2^2
\le |x'|^2 |z'|^2 - \< x'', z'' \> ^2
\le |z'|^2 - \< x'', z'' \> ^2
\end{equation}
because $|x'| \le \|x\|_2 \le 1$.

Now we apply \eqref{prime} and \eqref{double prime} to the running
approximation $x = x_{k-1}$ and the next approximation $x_1 = x_k$
obtained with a random $z = Z_k$.
Denoting $p_k = \langle x''_k, Z''_k \rangle$,
we have by \eqref{prime} that
$|x'_k - x'_{k-1}| \le |Z'_k|^2 + |p_k Z'_k|$
and by \eqref{double prime} that
$ \|x''_k\|_2^2 - \|x''_{k-1}\|_2^2 \le |Z'_k|^2 - |p_k|^2$.
Since $x'_0 = 1$ and $x''_0 = 0$, we have
\begin{equation}
|x'_k - 1|
\le \sum_{j=1}^k |x'_j - x'_{j-1}|
\le \sum_{j=1}^k |Z'_j|^2 + \sum_{j=1}^k |p_j Z'_j| \label{x'k}
\end{equation}
and
$$
\|x''_k\|_2^2
= \sum_{j=1}^k \big( \|x''_j\|_2^2 - \|x''_{j-1}\|_2^2 \big)
\le \sum_{j=1}^k |Z'_j|^2 - \sum_{j=1}^k |p_j|^2.
$$
Since $\|x''_k\|_2^2 \ge 0$, we conclude that
$\sum_{j=1}^k |p_j|^2
\le \sum_{j=1}^k |Z'_j|^2$.
Using this, we apply Cauchy-Schwartz inequality in \eqref{x'k}
to obtain
$$
|x'_k - 1|
\le \sum_{j=1}^k |Z'_j|^2
  + \Big( \sum_{j=1}^k |Z'_j|^2 \Big)^{1/2}
    \Big( \sum_{j=1}^k |Z'_j|^2 \Big)^{1/2}
= 2 \sum_{j=1}^k |Z'_j|^2.
$$
Since all $Z_j$ are copies of the random vector $Z$,
we conclude by \eqref{expectation Z'} that
$\E |x'_k - 1| \le 2k \, \E|Z'|^2 \le 2k/\k(A)^2$.
Thus
$\E \|x_k\| \ge \E|x'_k| \ge 1 - 2k/\k(A)^2$.
This proves the theorem, actually with the stronger conclusion
$$
\E \|x_k - x\|_2 \ge \big( 1 - 2k/\k(A)^2 \big) \cdot \|x_0 - x\|_2.
$$
The actual conclusion follows by Jensen's inequality.
\endproof

\subsection{The upper estimate is attained}

If $k(A) = 1$ (equivalently, if $\k(A) = \sqrt{n}$  by \eqref{condition numbers}),
then the estimate in Theorem~\ref{solution}
becomes an equality. This follows directly from the proof
of Theorem~\ref{solution}.

Furthermore, there exist arbitrarily large systems and with arbitrarily
large condition numbers $k(A)$ for which the estimate in Theorem~\ref{solution}
becomes an equality. Indeed, let $n$ and $m \ge n$ be arbitrary numbers. 
Let also $\k \ge \sqrt{n}$ be any number such that $m/\k^2$ is an integer.
Then there exists a system \eqref{system} of $m$ equations in $n$ variables
and with $\k(A) = \k$, for which the estimate in Theorem~\ref{solution}
becomes an equality for every $k$.

To see this, we define the matrix $A$ with the help of any orthogonal set
$e_1, \ldots, e_n$ in $\R^n$. Let the first $m/\k^2$ rows of $A$
be equal to $e_1$, the other rows of $A$ be equal to one of the vectors
$e_j$, $j > 1$, so that every vector from this set repeats at least
$m/\k^2$ times as a row (this is possible because $\k^2 \ge n$).
Then $\k(A) = \k$ (note that \eqref{rows of A} is attained for $z = e_1$).

Let us test our algorithm on the system $Ax = 0$ with the initial
approximation $x_0 = e_1$ to the solution $x = 0$. Every step of the algorithm
brings the running approximation to $0$ with probability $\k^{-2}$
(the probability of picking the row of $A$ equal to $e_1$ in uniform sampling),
and leaves the running approximation unchanged with probability $1 - \k^{-2}$.
By the independence, for all $k$ we have
$$
\E \|x_k - x_0\|_2^2 = \bigl(1 - \k^{-2} \bigr)^k \cdot \|x_0 - x\|_2^2.
$$

\section{Numerical experiments and comparisons}
\label{s:numeric}

\subsection{Reconstruction of bandlimited signals from nonuniform sampling} 
\label{ss:samp}

The reconstruction of a bandlimited function $f$ from its nonuniformly
spaced sampling values $\{f(t_k)\}$ is a classical problem in 
Fourier analysis, with a wide range of applications~\cite{BF01}.
We refer to~\cite{FG93,FGS95} for various efficient numerical algorithms.
Staying with the topic of this paper, we focus on the Kaczmarz method, 
also known as POCS (Projection Onto Convex Sets) method in signal 
processing~\cite{YS90}. 

As efficient finite-dimensional model, appropriate
for a numerical treatment of the nonuniform sampling problem,
we consider trigonometric polynomials~\cite{Gro99}. In this model the
problem can be formulated as follows:
Let $f(t) = \sum_{l=-r}^r x_l e^{2\pi i lt}$, where $x=\{x_l\}_{l=-r}^r \in
\Cst^{2r+1}$. Assume we are given the nonuniformly spaced nodes
$\{t_k\}_{k=1}^m$ and the sampling values $\{f(t_k)\}_{k=1}^m$.
Our goal is to recover $f$ (or equivalently $x$).

The solution space for the $j$-the equation is given by
the hyperplane 
$$\{y : \< y, D_r(\cdot -t_j) \> = f(t_j) \},$$
where $D_r$ denotes the Dirichlet kernel 
$D_r(t) = \sum_{k=-r}^r e^{2\pi i kt}$.
% \frac{1}{2M+1}\frac{\sin(2M+1)\pi t}{\sin(\pi t)}.$$
Feichtinger and Gr\"ochenig argued convincingly (see e.g.~\cite{FG93})
that instead of $D_r(\cdot - t_j)$ one should consider the weighted Dirichlet 
kernels $\sqrt{w_j} D_r(\cdot - t_j)$, where the weight 
$w_j = \frac{t_{j+1}-t_{j-1}}{2}, j=1, \dots, m$. 
The weights are supposed to compensate for varying density in the sampling
set.

Formulating the resulting conditions in the Fourier domain,
we arrive at the linear system of equations~\cite{Gro99}
\begin{equation}
Ax=b, \qquad \text{where $A_{j,k} = \sqrt{w_j} e^{2\pi i k t_j}$, 
\,\,\, $b_j = \sqrt{w_j} f(t_j)$},
\label{trigsys}
\end{equation}
with $j=1,\dots,m$; $k = -r,\dots,r$. 
Let use denote $n:=2r+1$ then $A$ is an $m \times n$ matrix.

%Formulating the resulting conditions in the Fourier domain,
%we arrive at the linear system of equations $Ax=b$, 
%where $A$ is an $m \times 2r+1$ matrix
%with entries $A_{j,k} = \sqrt{w_j} e^{2\pi i k t_j}$, $j=1,\dots,m,
%k = -M,\dots,M$, and $b$ is a vector of length $m$ with entries 
%$b_j = \sqrt{w_j} f(t_j)$. Let us denote $n:=2M+1$.

Applying the standard Kaczmarz method (the POCS method as proposed
in~\cite{YS90}) to~\eqref{trigsys} means that
we sweep through the projections in the natural order, i.e., 
we first project on the hyperplane associated with the first row
of $A$, then proceed to the second row, the third row, etc. 
As noted in~\cite{FG93} this is a rather inefficient way of implementing
the Kaczmarz method in the context of the nonuniform sampling problem.
It was suggested in~\cite{CFM92} that the convergence can be
improved by sweeping through the rows of $A$ in 
a random manner, but no proof of the (expected) rate of convergence
was given. \cite{CFM92} also proposed another variation of the Kaczmarz 
method in which one projects in each step onto that hyperplane 
that provides the largest decrease of the residual error. 
This strategy of {\em maximal correction} turned out to provide very good 
convergence, but was found to be impractical due to the enormous 
computational overhead, since in each step all $m$ projections have to be 
computed in order to be able to select the best hyperplane to project on. 
It was also observed in~\cite{CFM92} that this maximal correction
strategy tends to select the hyperplanes associated with large weights
more frequently than hyperplanes associated with small weights.

Equipped with the theory developed in Section~\ref{s:random} we can
shed light on the observations mentioned in the previous paragraph. 
Note that the $j$-th row of $A$ in~\eqref{trigsys} has squared norm equal 
to $n w_j$. Thus our Algorithm~\ref{alg} chooses the $j$-th row of $A$ with 
probability $w_j$. Hence Algorithm~\ref{alg} can be interpreted as
a probabilistic, computationally very efficient implementation of
the maximal correction method suggested in~\cite{CFM92}.

Moreover, we can give a bound on the expected rate of convergence
of the algorithm. Theorem~\ref{solution} states that this rate 
depends only on the scaled condition number $\kappa(A)$,
which is bounded by $k(A) \sqrt{n}$ by \eqref{condition numbers}.
The condition number $k(A)$ for the trigonometric system 
\eqref{trigsys} has been estimated by Gr\"ochenig \cite{Gro92}.
For instance we have the following
\begin{theorem}[Gr\"ochenig]
If the distance of every sampling point $t_j$ to its neighbor
on the unit torus is at most $\delta < \frac{1}{2r}$,
then $k(A) \le \frac{1+2\delta r}{1-2\delta r}$. In particular, if
$\delta \le \frac{1}{4r}$ then $k(A) \le 3$.
\end{theorem}
Furthermore we note that our algorithm can be straightforward applied 
to the approximation of multivariate trigonometric polynomials.
We refer to~\cite{BG} for condition number estimates for this case.

\medskip

\begin{figure}[ht!]
\begin{center}
\includegraphics[width=100mm,height=70mm]{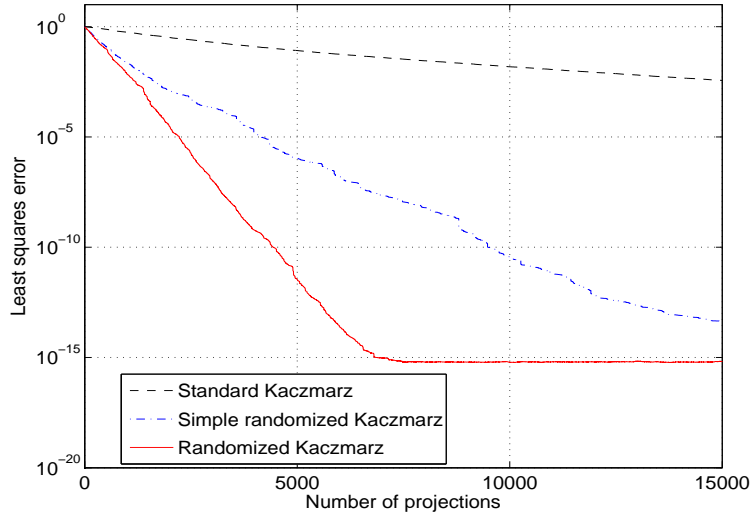}
\caption{Comparison of rate of convergence for the randomized Kaczmarz method
described in Algorithm~\ref{alg} and other Kaczmarz methods applied to
the nonuniform sampling problem described in the main text.}
\label{fig:trig}
\end{center}
\end{figure}

In our numerical simulation, we let $r=50, m = 700$ and generate the sampling
points $t_j$ by drawing them randomly from a uniform distribution in $[0,1]$
and ordering them by magnitude. 
We apply the standard Kaczmarz method, the randomized Kaczmarz method,
where the rows of $A$ are selected at random with equal probability
(labeled as {\em simple randomized Kaczmarz} in Figure~\ref{fig:trig}), and 
the randomized Kaczmarz method of Algorithm~\ref{alg} 
(where the rows of $A$ are selected at random
with probability proportional to the $2$-norm of the rows).
We plot the least squares error $\|x-x_k\|_2$ versus the number of
projections, cf.~Figure~\ref{fig:trig}. Clearly, Algorithm~\ref{alg}
significantly outperforms the other Kaczmarz methods, demonstrating
not only the power of choosing the projections at random, but also the
importance of choosing the projections according to their relevance.

\subsection{Comparison to conjugate gradient algorithm} \label{ss:CG}

In recent years conjugate gradient (CG) type methods have emerged as
the leading iterative algorithms for solving large linear systems of
equations, since they often exhibit remarkably fast 
convergence~\cite{GL96,Han95}.
How does Algorithm~\ref{alg} compare to the CG algorithms?

\if 0
It is not surprising that one can easily construct examples for which CG
(or its variations, such as CGLS or LSQR~\cite{Han95})
will clearly outperform the proposed method. For instance, take a matrix
whose singular values, all but one, are equal to one, while the remaining
singular value is very small, say $10^{-8}$. While this matrix is far 
from being well-conditioned, CGLS will nevertheless converge in
only two iterations, due to the clustering of the spectrum of $A$,
cf.~\cite{SV86}.
In comparison, the proposed Kaczmarz method will converge extremely slowly
in this example by  since $R \approx \eps^{-2}$ and thus 
$1-\frac{1}{R} \approx 1$.

On the other hand, Algorithm~\ref{alg} can outperform CGLS in other cases 
for which CGLS is actually quite well suited. As an example we consider the 
important case of a Gaussian random matrix $A$ of dimension $m \times n$ with 
$m \ge n$ and randomly generated $x$. We will now elaborate on this
example in more detail. 
\fi

The rate of convergence of CGLS applied to $Ax=b$ is bounded 
by~\cite{GL96}
\begin{equation}
\|x_k - x\|_{A^{\ast}A} \le 2 \|x_0 - x\|_{A^{\ast}A}
\Big(\frac{k(A)-1}{k(A)+1}\Big)^k,
\label{cgconv}
\end{equation}
where\footnote{Note that since we either need to apply CGLS to
$Ax=b$ or CG to $A^{\ast}Ax = A^{\ast}b$ we indeed have to 
use $k(A)=\sqrt{k(A^{\ast}A)}$ here and not $\sqrt{k(A)}$.
The asterisk $^*$ denotes complex transpose here.}
$\|y\|_{A^{\ast}A} := \sqrt{\langle Ay,Ay\rangle}$.

It is known that the CG method may converge faster when the 
singular values of $A$ are clustered~\cite{SV86}. For instance, take a matrix
whose singular values, all but one, are equal to one, while the remaining
singular value is very small, say $10^{-8}$. While this matrix is far 
from being well-conditioned, CGLS will nevertheless converge in
only two iterations, due to the clustering of the spectrum of $A$,
cf.~\cite{SV86}.
In comparison, the proposed Kaczmarz method will converge extremely slowly
in this example by Theorem~\ref{general lower estimate}, 
since $\k(A) \approx 10^{8}$.

On the other hand, Algorithm~\ref{alg} can outperform CGLS on problems
for which CGLS is actually quite well suited, in particular for 
random Gaussian matrices $A$, as we show below.

\paragraph{Solving random linear systems}
Let $A$ be a $m \times n$ matrix whose entries are 
independent $N(0,1)$ random variables.
Condition numbers of such matrices are well studied, when the 
aspect ratio $y := n/m < 1$ is fixed and the size $n$ of the matrix grows
to infinity. Then the following almost sure convergence was 
proved by Geman \cite{Gem80} and Silvestein \cite{Sil} respectively:
$$
\frac{\|A\|_2}{\sqrt{m}} \to 1 + \sqrt{y}; 
\qquad
\frac{1/\|A^{-1}\|_2}{\sqrt{m}} \to 1 - \sqrt{y}.
$$
Hence 
\begin{equation}				\label{k-random}
  k(A) \to \frac{1 + \sqrt{y}}{1 - \sqrt{y}}.
\end{equation}
Also, since $\frac{\|A\|_F}{\sqrt{mn}} \to 1$, we have
\begin{equation}				\label{kappa-random}
  \frac{\k(A)}{\sqrt{n}} \to \frac{1}{1 - \sqrt{y}}.
\end{equation}
For estimates that hold for each finite $n$ rather than in the
limit, see e.g. \cite{ES} and \cite{Ede92}.

Now we compare the expected computation complexities of the 
randomized Kaczmarz algorithm proposed in Algorithm~\ref{alg} and CGLS
to compute the solution within error $\eps$ for the 
system \eqref{system} with a random Gaussian matrix $A$.

We estimate the expected number of iterations (projections) $k_\e$ 
for Algorithm~\ref{alg} to achieve an accuracy $\eps$ 
in \eqref{k-eps}. Using bound \eqref{kappa-random}, 
we have 
$$
\E \, k_\e \approx \frac{2 n}{(1 - \sqrt{y})^2} \log \frac{1}{\e}
$$
as $n \to \infty$. 
Since each iteration (projection) requires $n$ operations, 
the total expected number of operations is 
\begin{equation}
  \text{Complexity of randomized Kaczmarz} \approx 
  \frac{2n^2}{(1-\sqrt{y})^2} \log \frac{1}{\e}.
\label{kaczcomp}
\end{equation}

The expected number of iterations $k'_\e$ for CGLS to 
achieve the accuracy $\eps$ can be estimated using \eqref{cgconv}.
First note that the norm $\|\cdot\|_{A^*A}$ is on average proportional 
to the Euclidean norm $\|z\|_2$. Indeed, for any fixed vector $z$ one has
$\E \|z\|_{A^*A}^2 = \E \|Az\|_2^2 = m \|z\|_2^2$. 
Thus, when using CGLS for a random matrix $A$, we can
expect that the bound \eqref{cgconv} on the convergence also 
holds for the Euclidean norm.

%THE LAST STATEMENT IS NOT RIGOROUS. THE VECTOR $x_k$ IS RANDOM: IT 
%DEPENDS ON $A$. THUS WE CAN NOT TAKE THE EXPECTATION IN \eqref{cgconv}
%AND USE THE ABOVE PROPORTIONALITY OF THE NORMS. BUT INTUITIVELY ITS OK.

Consequently, the expected number of iterations $k'_\e$ in CGLS
to compute the solution within accuracy $\e$ as in \eqref{within eps}
is 
$$
\E k_\e \approx \frac{\log \frac{2}{\e}}{\log K(A)}
\qquad
\text{where} \quad K(A) = \frac{k(A)+1}{k(A)-1}.
$$ 
By \eqref{k-random}, for random matrices $A$ of growing size we have 
$K \to 1/\sqrt{y}$ almost surely. Thus 
$$
\E k_\e \approx \frac{2\log \frac{2}{\e}}{\log \frac{1}{y}}.
$$
The main computational task in each iteration of CGLS consists of
two matrix vector multiplications, one with $A$ and one with $A^{\ast}$, 
each requiring $m \times n = n^2/y$ operations. 
Hence the total expected number of operations is 
\begin{equation}
  \text{Complexity of CGLS} 
  \approx \frac{4 n^2}{y \log \frac{1}{y}} \cdot \log \frac{2}{\e}.
\label{cgcomp}
\end{equation}

\begin{figure}[ht!]
\begin{center}
\includegraphics[width=100mm,height=70mm]{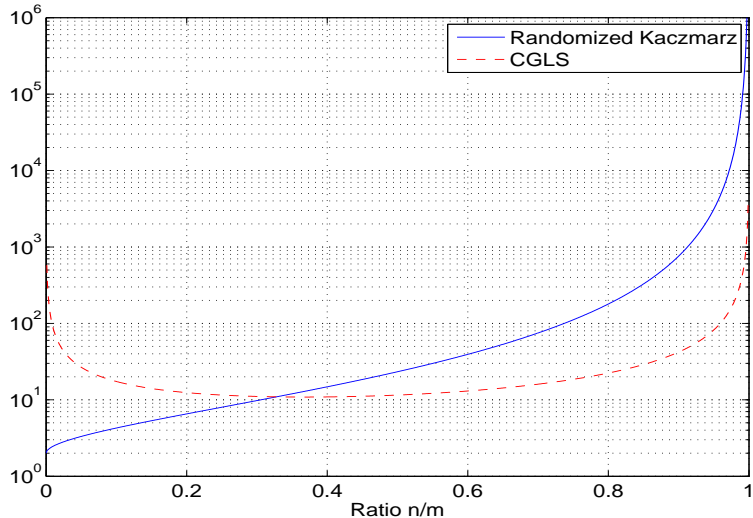}
\caption{Comparison of the computational complexities~\eqref{kaczcomp}
(randomized Kaczmarz method) and~\eqref{cgcomp} (conjugate gradient 
algorithm) as functions of the ratio $y=\frac{n}{m}$
(the common factors $n^2$ and $\log(1/\eps)$ in~\eqref{kaczcomp}
and~\eqref{cgcomp} are ignored in the two curves).}
\label{fig:cgkaczcomp}
\end{center}
\end{figure}

It is easy to compare the complexities~\eqref{kaczcomp} and~\eqref{cgcomp}
as functions of $y$, since $n^2$ and $\log(1/\eps)$ are common terms in both
(using the approximation $\log(2/\eps) \approx \log(1/\eps)$ for small $\eps$), 
cf.~also Figure~\ref{fig:cgkaczcomp}. A simple computation shows 
that \eqref{kaczcomp} and~\eqref{cgcomp} are essentially equal when 
$y\approx \frac{1}{3}$. Hence for Gaussian matrices 
our analysis predicts that Algorithm~\ref{alg} outperforms CGLS in terms 
of computational efficiency when $m > 3n$.

While the computational cost of Algorithm~\ref{alg} decreases as
$\frac{n}{m}$ decreases, this is not the case for CGLS. Therefore it
is natural to ask for the optimal ratio $\frac{n}{m}$ for CGLS 
for Gaussian matrices that minimizes its overall computational complexity.
It is easy to see that for given $\eps$ the 
expression in~\eqref{cgcomp} is minimized if $y=1/e$, where $e$ is
Euler's number. Thus if we are given an $m\times n$ Gaussian matrix
(with $m > en$), the most efficient strategy to employ CGLS is to first
select a random submatrix $\Ae$ of size $en \times n$ from $A$ (and 
the corresponding subvector $\be$ of $b$) and apply CGLS to the subsystem
$\Ae x= \be$. This will result in the optimal computational complexity 
$4e n^2 \log\frac{2}{\eps}$ for CGLS. 

Thus for a fair comparison between the randomized Kaczmarz method and CGLS,
we will apply CGLS in our numerical simulations
to both the ``full'' system $Ax=b$ as well as to a subsystem
$\Ae x = \be$, where $\Ae$ is an $en \times n$ submatrix of
$A$, randomly selected from $A$. 

%Nevertheless, our numerical experiments, presented below, reveal that 
% Algorithm~\ref{alg} performs actually even better than predicted.

\medskip
In the first simulation we 
let $A$ be of dimension $300 \times 100$, the entries of $x$ are
also drawn from a normal distribution. We apply both, CGLS and
Algorithm~\ref{alg}. We apply CGLS to the full system of size $300 \times
100$ as well as to a randomly selected subsystem of size $272 \times 100$ 
(representing the optimal size $e n \times n$, computed above).
Since we know the true solution we can compute the actual least squares
error $\|x - x_k\|$ after each iteration. Each method is terminated after
reaching the required accuracy $\eps = 10^{-14}$. We repeat the experiment 
100 times and for each method average the resulting least squares errors.

\begin{figure}[ht!]
\begin{center}
\includegraphics[width=100mm,height=70mm]{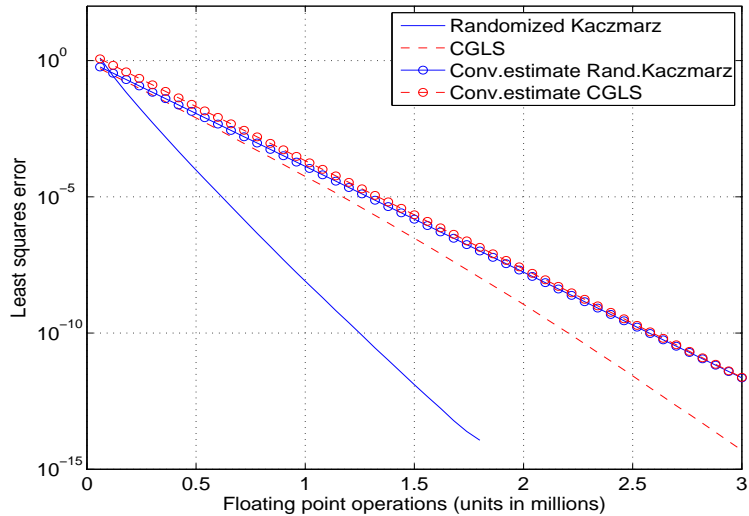}
\caption{Comparison of rate of convergence for the randomized Kaczmarz method
described in Algorithm~\ref{alg} and the conjugate gradient least
squares algorithm for a system of equations with a Gaussian 
matrix of size $300 \times 100$.}
\label{fig:cg1}
\end{center}
\end{figure}

In Figure~\ref{fig:cg1} we plot the averaged least squares error 
(y-axis) versus the number of floating point operations (x-axis), 
cf.~Figure~\ref{fig:cg1}. We also plot the estimated convergence rate for 
both methods. Recall that our estimates predict 
essentially identical bounds on the convergence behavior for CGLS and 
Algorithm~\ref{alg} for the chosen parameters ($m=3n$).
Since in this example the performance of CGLS applied to the full system 
of size $300\times 100$ is almost identical to that of CGLS applied to 
the subsystem of size $272 \times 100$, we display only the results
of CGLS applied to the original system.

While CGLS performs somewhat better than the (upper) bound predicts,
Algorithm~\ref{alg} shows a significantly faster convergence
rate. In fact, the randomized Kaczmarz method is almost twice as efficient
as CGLS in this example.

\begin{figure}[ht!]
\begin{center}
\includegraphics[width=100mm,height=70mm]{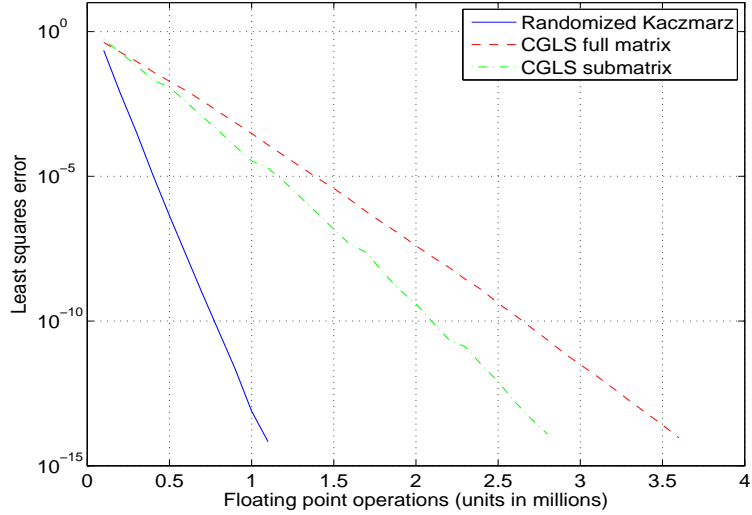}
\caption{Comparison of rate of convergence for the randomized Kaczmarz method
described in Algorithm~\ref{alg} and the conjugate gradient least
squares algorithm for a system of equations with a Gaussian 
matrix of size $500 \times 100$.}
\label{fig:cg2}
\end{center}
\end{figure}

In the second example we let $m=500, n=100$. In the same way as before,
we illustrate the convergence behavior of CGLS and Algorithm~\ref{alg}.
In this example we display the convergence rate for CGLS applied to
the full system (labeled as CGLS full matrix) of size $500\times 100$ as well 
as to a random subsystem of size $272 \times 100$ (labeled as CGLS
submatrix). As is clearly visible in~Figure~\ref{fig:cg2} CGLS applied to 
the subsystem provides better performance than CGLS applied to the
full system, confirming our theoretical analysis. Yet again, 
Algorithm~\ref{alg} is even more efficient
than predicted, this time outperforming CGLS by a factor of 3 (instead of
a factor of about 2 according to our theoretical analysis).

\noindent
{\bf Remark:} An important feature of the conjugate gradient algorithm
is that its computational complexity reduces significantly 
when the complexity of the matrix-vector multiplication is
much smaller than $\ord(mn)$, as is the case e.g.\ for
Toeplitz-type matrices. In such cases conjugate gradient algorithms
will outperform Kaczmarz type methods.

\section{Some open problems} \label{s:open}

In this final section we briefly discuss a few loose ends and some open
problems.

\medskip
\noindent
{\bf Kaczmarz method with relaxation:} It has been observed that
the convergence of the Kaczmarz method can be accelerated by introducing
relaxation.
In this case the iteration rule becomes
\begin{equation}
\label{kcazrelax}
x_{k+1}=x_k +\lambda_{k,i} \frac{b_i - \langle a_i,x_k\rangle}{\|a_i\|_2^2}a_i,
\end{equation}
where the $\lambda_{k,i}$, $i=1,\dots,m$ are relaxation parameters. 
For consistent systems the relaxation parameters must satisfy~\cite{HLL78}
\begin{equation}
0 < \liminf_{k \to \infty} \lambda_{k,i}\le \limsup_{k \to \infty}
\lambda_{k,i} < 2
\label{relax}
\end{equation}
to ensure convergence. 

\begin{figure}[ht!]
\begin{center}
\includegraphics[width=100mm,height=70mm]{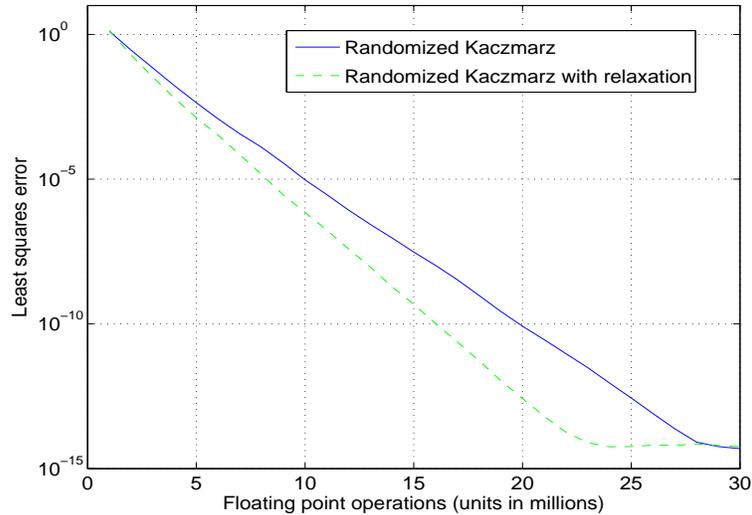}
\caption{Comparison of rate of convergence for the randomized Kaczmarz method
with and without relaxation parameter. We have used 
$\lambda = 1+\frac{n}{m}$ as relaxation parameter.}
\label{fig:relax}
\end{center}
\end{figure}

We have observed in our numerical simulations that for instance for
Gaussian matrices a good choice for the relaxation parameter is to
set $\lambda_{k,i} :=\lambda = 1+\frac{n}{m}$ for all $k$ and $i$.
While we do not have
a proof for an improvement of performance or even optimality,
we provide the result of a numerical simulation that is typical
for the behavior we have observed, cf.\ Figure~\ref{fig:relax}.

\medskip
\noindent
{\bf Inconsistent systems:} Many systems arising in practice are
inconsistent due to noise that contaminates the right hand side. In this 
case it has been shown that convergence to the least squares solution can 
be obtained with (strong under)relaxation~\cite{CEG83,HN90}.
We refer to~\cite{HN90,HN90a} for suggestions for the choice of
the relaxation parameter as well as further in-depth analysis for this case.

While our theoretical analysis presented in this paper assumes consistency 
of the system of equations, it seems quite plausible that the randomized 
Kaczmarz method combined with appropriate underrelaxation should 
also be useful for inconsistent systems.

\if 0 
\bibliographystyle{plain}
\bibliography{nuhag,mathbook,thomas,linalg,samp,random,general}
\fi

\end{document}